%% file: CHARREV.TEX
\input mydef

\def\cite#1{[#1]}
\def\frac#1#2{{#1\over #2}}
\def\Kasten{\hfill\hbox{             
    \vrule width 0.3pt height 4.7pt depth 0.3pt
    \vrule width 4.4pt height 0pt   depth 0.3pt \hskip -4.4pt
    \vrule width 4.4pt height 4.7pt depth -4.4pt
    \vrule width 0.3pt height 4.7pt depth 0.3pt \hskip 2pt}}

\font\smspec=cmcsc10 scaled \magstep0
\rightline{\it February 1994}
\bn\bn\bn
\cline{\spec Generalized Functionals in Gaussian Spaces}
\bn
\cline{\smspec -- The Characterization Theorem Revisited --}
\bn\bn
\cline{\it Yu.G. Kondratiev$^1$, P. Leukert$^1$, J. Potthoff$^3$, 
L. Streit$^{1,2}$, W. Westerkamp$^1$}
\bn
$$\vbox{\hbox {$\,^1$ {\smspec BiBoS}, Universit\"at Bielefeld}
       \hbox{$\,^2$ Centro das Ci\^encias Matem\'aticas, 
             Universidade da Madeira}
       \hbox  {$\,^3$ Lehrstuhl f\"ur Mathematik V, Universit\"at Mannheim}}
$$
\bn
\cline{Published in: Journal of Functional Analysis 141 No 2 (1996)}
\bn
\Th{Abstract.}
Gel'fand triples of test and generalized functionals in Gaussian spaces are
constructed and characterized.

\section{1. Introduction}
\noindent
In recent years there was an increasing interest in white noise analysis,
due to its rapid developments in mathematical structure and applications in
various domains. Especially, the circle of ideas going under the heading
`characterization theorems' has played quite an important role in the last
few years. These results \cite{21}, \cite{33}, \cite{42}, and their
variations and refinements (see, e.g., \cite{31}, \cite{34}, \cite{37}, 
\cite{48}, \cite{50}, \cite{54}, and references quoted there), provide
a deep insight into the structure of spaces of smooth and generalized random
variables over the white noise space or -- more generally -- Gaussian spaces.
Also, they allow for rather straightforward applications of these notions to
a number of fields: for example, Feynman integration \cite{11}, \cite{16},
\cite{20}, \cite{32}, representation of quantum field theory \cite{2},
\cite{43}, stochastic equations \cite{7}, \cite{30}, \cite{39}, 
\cite{40}, \cite{41}, intersection local times \cite{10}, \cite{49},
Dirichlet forms \cite{3}, \cite{4}, \cite{15}, infinite
dimensional harmonic analysis \cite{14} and so forth. Moreover,
characterization theorems have been at the basis of new methods for the
construction of smooth and generalized random variables \cite{24}, \cite
{34} which seem to be useful in applications untractable by existing
methods (e.g., \cite{18}, \cite{19}).
\medskip
The purpose of the present article is four-fould: We wish 1.\ to clarify and
generalize the structure of the existing characterization theorems, and at
the same time, 2.\ to review and unify recent developments in this direction,
3.\ to establish the connection to rich, related mathematical literature 
[1], \cite{8}, \cite{12}, [45], [53],  which might be helpful in 
future developments, and -- last but not least -- 4.\ to fill a gap in the 
article \cite{42}. In this
sense, the present paper attempts to give known results and some existing
`folklore' around them a general form that can be used as a reference for
future research. In the course of doing this, we also establish some new
results, for instance an analytic extension property of $U$--functionals,
and the topological invariance of certain spaces of
generalized random variables with respect to different construction schemes.
\medskip
This article is organized as follows. In Section 2 we present some
notions and results from complex analysis on topological vector
spaces, from Gaussian analysis, and from Fock space theory. In
particular, we construct a nuclear rigging
$$
(\NC)\,\subset\,\G(\HC)\,\subset\,(\NC)^*
$$ 
of the symmetric Fock space $\G(\HC)$ over a Hilbert space $\HC$.
We give a construction of the second quantized space $({\cal N})$ solely in
terms of the topology of ${\cal N}$, independent of the particular
representation as a projective limit. Via the well-known 
Wiener--It\^o--Segal isomorphism this provides a rigging of the $L^2$--space
over a Gaussian measure space by spaces of smooth and generalized
random variables. In Section 3 we study the  $U$--functionals associated
with the elements in $(\NC)$ and $(\NC)^*$. We derive an analytic
extension property for $U$--functionals, and use this to prove theorems 
which characterize $(\NC)$ and $(\NC)^*$ in terms of their 
$\SC$--transforms. In Section 4 we prove two corollaries of the
characterization theorem for $(\NC)^*$ which appear to be useful in 
applications.

\section{2. Preliminaries}
\noindent
{\bf 2.1 G--Entire Functions}
\sn
We provide some well-known facts from complex analysis on topological vector
spaces (see, e.g., \cite{8} and \cite{12}) with a view towards
applications in the next section.
\medskip
Let ${\cal E}$ be a locally convex complex vector space.

\Th{Definition 1.}{\it 
A mapping P from ${\cal E}$ into ${\C}$ is called an {\bf n--homogeneous
polynomial}, if it is the composition of the diagonal mapping $\Delta
_n:x\mapsto (x,x,\ldots, x)$ from $\EC$ into $\EC^n$ and a symmetric 
$n$--linear mapping $L$ from $\EC^n$ into $\C$, i.e., 
$P=L\circ \Delta _n\equiv \widehat L$. Let $P_n({\cal E})$ denote the space
of all $n$--homogeneous polynomials.} 

\Th{Definition 2.}{\it
A function F defined on ${\cal E}$ with values in 
${\C}$ is said to be {\bf G--entire} if for all $\xi, \eta\in 
{\cal E}$ the complex valued function 
$$
z\mapsto F(\eta+z\xi),\quad z\in\C, 
$$
is entire. Let $H_G({\cal E})$ denote the set of all G--entire mappings from $%
{\cal E}$ into ${\C}$. 
}
\medskip
For $F\in H_G(\EC),\,\eta\in\EC$, there exists a unique sequence
$({1\over n!}\wh{d^nF(\eta)},\,n\in\N_0)$ of homogeneous polynomials
${1\over n!}\wh{d^nF(\eta)}\in P_n(\EC),\,n\in\N$, such that for all
$\xi\in\EC$,
$$
F(\eta+\xi)= \sum_{n=0}^\infty {1\over n!}\,\wh{d^n F(\eta)}(\xi).
\eqno{(1)}
$$
Of course, $\widehat{d^nF(\eta)}( \xi)$ is the
$n$--th partial derivative of $F$ at $\eta$ in the direction $\xi$. The
corresponding $n$--linear form is denoted by 
$d^nF(\eta)(\xi_1,\ldots,\xi_n),\,\xi_1,\ldots,\xi_n\in\EC$.

\Th{Definition 3.}{\it 
Let $F$ be a mapping from $\EC$ into $\C$. $F$ is called {\bf entire}
if it is in $H_G(\EC)$, and if it is continuous. $H({\cal E})$
denotes the space of all entire functions on $\EC$.
}
\medskip\noindent

\mn
{\bf Proposition 4.}\quad{\it 
Let $F\in H_G({\cal E})$. Then $F\in H({\cal E})$ if and
only if $F$ is locally bounded.
}
\medskip
We conclude this subsection by stating a result which is related to the
celebrated ``cross theorem" of Bernstein. For a review of such results
we refer the interested reader also to [1]. The following is
a special case of a result by  Siciak: if we make use of the  
fact that any segment of the real line in the complex plane has strictly
positive transfinite diameter, then Corollary 7.3 in [45]
implies 

\Th{Proposition 5.}{\it 
Let $n\in\N,\,n\ge 2$, and $f$ be a complex valued function on $\R^n$.
Assume that for all $k=1,2,\ldots,n$, and 
$(x_1,\ldots,x_{k-1},x_{k+1},\ldots,x_n)\in\R^{n-1}$, the mapping
$$
x_k\lms f(x_1,\ldots,x_{k-1},x_k,x_{k+1},\ldots,x_n),
$$
from $\R$ into $\C$ has an entire extension. Then $f$ has an entire
extension to $\C^n$.
}

\bn
{\bf 2.2 Gaussian Spaces}
\sn
The primordial object of Gaussian analysis (e.g., [6], \cite{17}, \cite{21}, 
\cite{22}, \cite{25}, \cite{26}, \cite{27}, \cite{28}, \cite{29}) is
a real separable Hilbert space ${\cal H}$. One then considers a 
rigging of 
${\cal H}$, ${\cal N}\subset {\cal H}\subset 
{\cal N}^*$, where ${\cal N}$ is a real nuclear
space (see below and \cite{13}), densely and continuously embedded into
$\HC$, and ${\cal N}^*$ is its dual ($\HC$ being identified with
its dual). A typical example (which appears for instance in white noise
analysis) is the rigging $\sr\subset L^2(\R)\subset\ssr$ of $L^2(\R)$
(with Lebesgue measure) by the Schwartz spaces of test functions and 
tempered distributions. 
\medskip
Via Minlos' theorem the canonical
Gaussian measure $\mu $ on ${\cal N}^{*}$ is introduced by giving
its characteristic function
$$
C(f)=\int_{{\cal N}^*}e^{i\left\langle \omega
,f\right\rangle }\,d\mu (\omega )=
e^{-\frac 12\left| f\right| _{\cal H}^2},\quad f\in {\cal N}. 
$$
The space $L^2({\cal N}^*,d\mu )\equiv (L^2)$ of 
(equivalence classes of) complex valued functions on 
${\cal N}^ *$ which are square-integrable with respect to 
$\mu$ has the well-known Wiener--It\^o--Segal chaos decomposition \cite{36}, 
\cite{46}, \cite{47}, and one has the familiar Segal isomorphism ${\cal I}$
between $(L^2)$ and the complex Fock space $\Gamma ({\cal H)}$ over the
complexification $\HC_\Ckl$ of $\HC$.
\medskip
Spaces of smooth functions on $\NC^*$ can be constructed by mapping 
appropriate subspaces of 
$\Gamma ({\cal H)}$ into $(L^2)$ via the unitary mapping 
${\cal I}^{-1}:\Gamma ({\cal H)\rightarrow }(L^2)$, see, e.g., the 
construction using
second quantized operators in [6], \cite{17}. In the present context, we prefer
to work exclusively in Fock space. At all times a `translation' into
function space language via ${\cal I}^{-1}$ is of course equally valid. 
\bn\goodbreak\noindent
{\bf 2.3 Second Quantized Spaces}
\sn
Our starting point is a real separable nuclear space $\NC$. It is well-known
(e.g., [38]) that the topology of $\NC$ is equivalent to the projective
limit topology of an increasing countable system $(|\,\cdot\,|_p,\,p\in
\N_0)$ of compatible Hilbertian norms $|\,\cdot\,|_p$. In other words,
$\NC$ is a countably Hilbert space [13],
$$
{\cal N}\equiv \bigcap_p{\cal H}_p, 
$$
where $\HC_p$ is equal to the completion of $\NC$ with respect to
$|\,\cdot\,|_p$. Moreover, the usual Hilbert--Schmidt property for the
embeddings holds, i.e., for every $p\in\N_0$ there exists a $p^\prime
>p$ so that the embedding $\iota$ from $\HC_{p^\prime}$ into
$\HC_p$ is a Hilbert--Schmidt operator. We shall denote the bilinear
dual pairing on ${\cal N}^{*}\times {\cal N}$ by 
$\left\langle \cdot ,\cdot
\right\rangle$. A special role -- also with a view towards Gaussian spaces --
is played by the Hilbert space ${\cal H}_0$, which we also denote by
${\cal H}$.
\medskip
Consider the Fock space, e.g., \cite{9},
[46], over ${\cal H}$%
$$
\Gamma \left( {\cal H}\right) =
\bigoplus_{n=0}^\infty \,{\cal H}_\Ckl^{\widehat{\otimes }n}, 
$$
where ${\cal H}_\Ckl^{\widehat{\otimes }n}$ is the symmetric $n$--fold
tensor product of ${\cal H}_\Ckl$ with itself. The 
inner product and norm $\left\|\, \cdot\, \right\| _0$ of $\Gamma({\cal H})$
are generated by
$$
\left( \varphi ^{\ten n},\psi ^{\ten n}\right) _{\Gamma \left( 
{\cal H}\right) }=n!\left( \varphi ,\psi \right) _{{\cal H}_\Ckkl}^n. 
$$
Likewise, for $p\in\N_0$, $n\in{\N}$, 
${\cal H}^{\widehat\otimes n}_{\Ckl,p}$ denotes the $n$--fold symmetric tensor 
product of ${\cal H}_{\Ckl,p}$ with itself, and it is considered as a subspace
of ${\cal H}^{\otimes n}_{\Ckl,p}$. The canonical norm of the latter is denoted
by $|\,\cdot\,|_p$, too (the meaning will be clear from the context).
The duals of ${\cal H}^{\otimes n}_{\Ckl,p}$ and 
${\cal H}^{\widehat\otimes n}_{\Ckl,p}$,
respectively, are denoted by ${\cal H}^{\otimes n}_{\Ckl,-p}$ and
${\cal H}^{\widehat\otimes n}_{\Ckl,-p}$, respectively. The Hilbertian norm of
${\cal H}^{\otimes n}_{\Ckl,-p}$ is denoted by $|\,\cdot\,|_{-p}$, and we remark
that for $n\in {\N},\,p\in {\Z},\, 
\Phi^{(n)}\in {\cal H}_{\Ckl,p}^{\otimes n}$, $|\Phi^{(n)}|_p$ is equal to the 
Hilbert--Schmidt norm of $\Phi^{(n)}$ considered as a linear form on
${\cal H}_{\Ckl,-p}^{\otimes n}$.  
For $q\in\N_0$, we introduce Hilbert spaces $\Gamma
_q\left( {\cal H}_p\right) $ as the completions of the space of finite 
direct sums%
$$
\textstyle{\bigoplus^\prime_n}\,{\cal H}_{\Ckl,p}^{\widehat{\otimes }n} 
$$
with respect to the inner product determined by 
$$
\left( \varphi ^{\ten n},\psi ^{\ten n}\right)
_{\Gamma _q\left( {\cal H}_p\right) }=2^{nq}n!\left( \varphi ,\psi \right)_{%
{\cal H}_{\Ckkl,p}}^n, 
\eqno{(2)}
$$
and denote the corresponding norms by $\left\|\, \cdot\, \right\| _{p,q}$. 
Finally we set%
$$
\left( {\cal N}\right) =\bigcap_{p,q}\,\Gamma _q\left( {\cal H}_p\right), 
$$
equipped with the projective limit topology.
\medskip
\noindent {\sl Remarks.} Evidently substitution of the value 2 in 
equation (2)
 by any other number strictly larger than 1 produces the same space $%
\left( {\cal N}\right)$. We use the same notation $\left( {\cal N}\right)$ 
for the nuclear subspace of $(L^2)$ corresponding to $%
\left( {\cal N}\right) $ under the Wiener--It\^o--Segal isomorphism.

\Th{Lemma 6.}{\it
$({\cal N})$ is nuclear.
}
\mn
{\sl Proof.}\quad Nuclearity of $\left( {\cal N}\right) $ follows
essentially from that of ${\cal N}.$ For fixed $p,q$ consider the embedding 
$$
I:\Gamma _{q^{\prime }}\left( {\cal H}_{p^{\prime }}\right) \rightarrow
\Gamma _q\left( {\cal H}_p\right) 
$$
where $p^{\prime }$ is chosen such that the embedding%
$$
\iota:\,{\cal H}_{p^{\prime }}\lra {\cal H}_p 
$$
is Hilbert--Schmidt. Then%
$$
I=\bigoplus_n \iota^{\ten n}. 
$$
Its Hilbert--Schmidt norm is easily estimated by using an orthonormal
basis (cf., e.g., \cite{17},
Appendix A.2) of $\Gamma_{q^\prime}({\cal H}_{p^\prime})$. The result is 
the bound
$$
\left\| I\right\| _{\rm HS}^2\le\sum_{n=0}^\infty 2^{n(q-q^{\prime })}
\left\|\iota\right\| _{\rm HS}^{2n} 
$$
which is finite for suitably chosen $q^{\prime }$.\Kasten
\medskip

\Th{Theorem 7.}{\it
The topology on $\left( {\cal N}\right) $ is uniquely determined by the
topology on ${\cal N}$.
}
\mn
{\sl Proof.}\quad Let us assume that we are given two different
systems of Hilbertian norms $\left|\, \cdot\, \right| _p$ and 
$\left|\, \cdot\, \right|_k^\prime$ , such that they induce the same 
topology on ${\cal N}$ .
For fixed $k$ and $l$ we have to estimate 
$\left\|\, \cdot\, \right\|_{k,l}^\prime$ by 
$\left\|\, \cdot\, \right\| _{p,q}$ for some $p,q$ (and
vice versa which is completely analogous). Since 
$\left|\, \cdot\, \right|_k^\prime$ has to be continuous 
with respect to the projective limit
topology on ${\cal N}$, there exists $p$ and a constant $C$ such that 
$\left|f\right| _k^\prime\leq C\left| f \right| _p$, for
all $f\in {\cal N}$, i.e., the injection $\iota$ from 
${\cal H}_{p}$ into the completion ${\cal K}_{k}$ of ${\cal N}$ with
respect to $|\,\cdot\,|_k^\prime$ is a mapping bounded by $C$.
We denote by $\iota$ also its linear extension from
$\HC_{\Ckl,p}$ into $\KC_{\Ckl,k}$. It follows from a straightforward 
modification of the  proof
of the Proposition on p.\ 299 in \cite{44}, that $\iota^{\otimes n}$
is bounded by $C^n$ from ${\cal H}_{\Ckl,p}^{\otimes n}$ into
${\cal K}_{\Ckl,k}^{\otimes n}$. Now we
choose $q$ such that $2^{\frac{q-l}2}\geq C$. Then
$$
\eqalign{
\left\|\, \cdot\, \right\| _{k,l}^{\prime 2}&=\sum_{n=0}^\infty n!\,
2^{nl}\left|\, \cdot\, \right|_k^{\prime 2}\cr
&\leq \sum_{n=0}^\infty n!\,2^{nl}C^{2n}\left|\, \cdot\, \right| _p^2\cr
&\leq \left\|\, \cdot\, \right\| _{p,q}^2,\cr 
}
$$
which had to be proved.\Kasten
\medskip
From general duality theory on nuclear spaces we know that the dual of $%
\left( {\cal N}\right) $ is given by%
$$
\left( {\cal N}\right) ^{*}=\bigcup_{p,q}\Gamma _q\left( {\cal H}_p\right)
^{*} 
$$
and one verifies that%
$$
\Gamma _q\left( {\cal H}_p\right) ^{*}=\Gamma _{-q}\left( {\cal H}%
_{-p}\right) . 
$$
We shall denote the bilinear dual pairing on $({\cal N})^{*}\times ({\cal N}%
) $ by 
$\left\langle\! \left\langle \cdot ,\cdot \right\rangle\! \right\rangle :$%
$$
\left\langle\! \left\langle \Phi ,\varphi \right\rangle\!\right\rangle
=\sum_{n=0}^\infty n!\,\langle \Phi ^{(n)},\varphi ^{(n)}\rangle, 
$$
where $\Phi\in\Gamma_{-q}({\cal H}_{-p})$ corresponds to the sequence 
$(\Phi^{(n)},\,n\in {\N}_0)$ with $\Phi^{(0)}\in {\C}$, and 
$\Phi^{(n)}\in {\cal H}^{\widehat\otimes n}_{\Ckl,-p},\,n\in{\N}$.
\medskip
\noindent {\sl Remark.} Consider the particular choice $\NC=\sr$.
Then $(\NC)^{(*)}$ coincide with (the Fock space equivalents of)
the well-known spaces $(\SC)^{(*)}$ of white noise functionals, see, e.g., 
[17], [42]. For the norms
$\|\vp\|_p\equiv \|\G(A^p)\vp\|_0$ introduced there, we have 
$\|\,\cdot\,\|_p=\|\,\cdot\,\|_{p,0}$, and $\|\,\cdot\,\|_{p,q}\le
\|\,\cdot\,\|_{p+{q\over 2}}$. More generally, if the norms on
$\NC$ satisfy the additional assumption that for all $p\ge 0$
and all $\e>0$ there exists $p'\ge 0$ such that $|\,\cdot\,|_p
\le\e|\,\cdot\,|_{p'}$, then the construction of Kubo and Takenaka
[26] (and other authors) leads to the same space $(\NC)$. The
construction presented here has the advantage of being manifestly
independent of the choice of any concrete system of Hilbertian
norms topologizing $\NC$.
\medskip
For the exponential vectors%
$$
\phi _f:= \sum_{n=0}^\infty \frac 1{n!}\,f^{\ten n} 
$$
one calculates the norms%
$$
\left\| \phi _f\right\| _{p,q}^2=e^{2^q\left| f\right| _p^2}, 
$$
and hence for all $f\in\NC$ they are in $\left( {\cal N}\right)$. 
This then allows for the following
\font\bi=cmbxti10 scaled \magstep0

\Th{Definition 8.}{\it
Let $\Phi\in ({\cal N})^*$. The {\bf {\bi S}--transform of $\Phi$} is the
mapping from $\NC$ into $\C$ given by
$$
S\Phi(f):=\langle\!\langle\Phi,\phi_f\rangle\!\rangle,\quad f\in\NC.
$$
}
\mn
 We note that the  exponential vectors $\{\phi_f,\,f\in\NC\}$, are a total 
set ${\cal E}$ in 
$({\cal N})$, and hence elements of $\left( {\cal N}\right) ^{*}$ are 
characterized by  their $S$--transforms. Furthermore, it is obvious 
that the $S$--transform of $\Phi\in(\NC)^*$ extends
to $\NC_\Ckl$: for $\xi\in\NC_\Ckl$ set $S\Phi(\xi)=
\langle\!\langle\Phi,\phi_\xi\rangle\!\rangle$, where $\phi_\xi$
is the complex exponential vector $\sum_n {1\over n!} \xi^{\ten n}\in
(\NC)$. 

\section{3. {\bi U}--Functionals and the Characterization Theorems}
\noindent
We begin with a definition.

\Th{Definition 9.}{\it
Let $F:{\cal N}\lra {\C}$ be such that 
\medskip
\iitem{C.1} for all $f,\,g\in\NC$, the mapping $\l\lms F(g+\l f)$
from $\R$ into $\C$ has an entire extension to $\l\in\C$,
\iitem{C.2} for some continuous quadratic form $B$ on $\NC$ there
exists constants $C,\,K>0$ such that for all $f\in\NC,\,z\in\C$,\par
$$
|F(zf)|\le C\,\exp(K\,|z|^2|B(f)|).
$$
Then F is called a {\bf {\bi U}--functional}.
}
\mn
{\sl Remark.}\quad Condition C.2 is actually equivalent to the more
conventional
\medskip{\it
\iitem{C.2$\,^\prime$} there exists constants $C,\,K>0$ and $p\in\N_0$,
so that for all $f\in\NC,\,z\in\C$,\par
$$
|F(zf)|\le C\,\exp(K\,|z|^2|f|^2_p).
\eqno{(3)}
$$
}
\Th{Lemma 10.} {\it
Every $U$--functional $F$ has a unique extension to an entire function
on $\NC_\Ckl$. Moreover, if the bound on $F$ holds in the form (3) then  
for all $\rho\in(0,1)$,
$$
|F(\xi)|\le C'\exp(K^\prime\,|\xi|^2_p),\quad\xi\in\NC_\Ckl,
$$
with $C'=C(1-\rho)^{-\half},\,K'=2\rho^{-1}e^2K$.
}
\mn
{\sl Proof.}\quad First we show that a $U$--functional $F$ has a G--entire
extension. The extension of $F$ (denoted by the same symbol)
is given by $F(\eta)=F(g_0+zg_1),\,\eta=g_0+zg_1\in\NC_\Ckl,\,g_0,g_1
\in\NC,\,z\in\C$. Let $\xi\in\NC_\Ckl$ be of the form $\xi=g_2+ig_3,\,
g_2,g_3\in\NC$. Consider the mapping
$$
(\l_1,\l_2,\l_3)\,\lms\,F(g_0+\l_1g_1+\l_2g_2+\l_3 g_3),
$$
from 
$\R^3$ into $\C$. Condition C.1 and Proposition 5 imply that this
function has an entire extension to $\C^3$. In particular, $F$ is
G--entire on $\NC_\Ckl$.
\medskip
Let $\xi\in\NC_\Ckl$, and consider the Taylor expansion of
$F(\xi)$ at the origin (cf.\ (1)):
$$
F(\xi)=\sum_{n=0}^\infty {1\over n!}\, \wh{d^nF(0)}(\xi).
\eqno{(4)}
$$
For all $f\in\NC,\,n\in\N,\,R>0$, we obtain from C.2$^\prime$ and
Cauchy's inequality the estimate
$$
|\wh{d^nF(0)}(f)|\le C\,n!\,R^{-n} e^{R^2K|f|_p^2}.
$$
We choose $R=({n\over 2K})^\half$, and get for $f\in\NC$ with $|f|_p=1$ the
inequality
$$
|\wh{d^nF(0)}(f)|\le C\,n!\,\Big({2eK\over n}\Big)^{n/2}.
$$
A standard polarization argument (see, e.g., [35, \S 3])  and homogeneity
of $\wh{d^nF(0)}$ yield the following bound for the $n$--linear
form $d^nF(0)$:
$$
|d^nF(0)(f_1,\ldots,f_n)|\le C(n!\,(2e^2K)^n)^{\half}\prod_{k=1}^n
                |f_k|_p,
\eqno{(5)}
$$
where $f_1,\ldots,f_n\in\NC$ (and we used ${n^n\over n!}\le e^n$).
\medskip
Since $d^nF(0)$ is $n$--linear on $\NC_\Ckl$, the last inequality
gives the estimate
$$
|d^nF(0)(\xi_1,\ldots,\xi_n)|\le C\,(n!\,(4e^2K)^n)^\half
                \prod_{k=1}^n|\xi_k|_p,
                \eqno{(6)}
$$
for $\xi_1,\ldots,\xi_n\in\NC_\Ckl$. In particular, the Taylor
coefficients in (4) have absolute value bounded by 
$$
C\,\Big({(4e^2K|\xi|^2_p)^n\over n!}\Big)^\half,
$$
and we get (by Schwarz' inequality) the following estimate for all
$\rho\in(0,1)$,
$$
|F(\xi)|\le  C(1-\rho)^{-\half}e^{2\rho^{-1}e^2K|\xi|_p^2},
\quad \xi\in\NC_\Ckl.
$$
Hence $F$ is locally bounded on $\NC_\Ckl$, and therefore
Proposition 4 implies that $F$ is entire.\Kasten
\medskip
Now we are ready to prove the following generalization of the main
result in [42] which characterizes the space $(\NC)^*$ in terms
of its $S$--transform.

\Th{Theorem 11.}{\it 
 A mapping $F:{\cal N}\rightarrow \C$ is the $S$%
--transform of an element in $({\cal N})^{*}$ if and only if it
is a U--functional.}
\mn
{\sl Proof.}\quad Let $\Phi\in(\NC)^*$. Then $\Phi\in\G_{-q}(\HC_{-p})$ for some
$p,\,q\in\N_0$. As we have remarked at the end of Section 2,
the $S$--transform of $\Phi$ extends to $\NC_\Ckl$, and therefore
it makes sense to consider the mapping $\xi\mapsto S\Phi(\xi)$
from $\NC_\Ckl$ into $\C$. We shall show that this mapping
is entire. We have
$$
S\Phi(\xi) = \sum_{n=0}^\infty\langle\Phi\nn,\xi^{\ten n}\rangle,
\quad\xi\in\NC_\Ckl.
$$
We estimate as follows:
$$
\eqalign{|S\Phi(\xi)|&\le \sum_{n=0}^\infty |\Phi\nn|_{-p}|\xi|_p^n\cr
        &\le\Big(\sum_{n=0}^\infty n!\,2^{-qn}|\Phi\nn|^2_{-q}\Big)^\half
         \Big(\sum_{n=0}^\infty {1\over n!} 2^{qn} |\xi|^{2n}_p\Big)^\half\cr
         &=\|\Phi\|_{-p,-q}\, e^{2^{q-1}|\xi|_p^2}.\cr
         }
$$
The last estimation shows that the power series for $S\Phi$
on $\NC_\Ckl$ converges uniformly on every bounded neighborhood of
zero in $\NC_\Ckl$, and therefore it defines an entire function
on this space [12]. In particular, C.1 holds for $S\Phi$. Moreover,
the choice $\xi=zf,\,z\in\C,\,f\in\NC$, shows that also C.2$^\prime$
is fulfilled. Hence $S\Phi$ is a $U$--functional.         

\medskip                                
Conversely let $F$ be a $U$--functional. We may assume the bound in the
form (3). Consider the $n$--linear form $d^nF(0)$ on $\NC_\Ckl$
constructed in the proof of Lemma 10. The estimate (6) shows that
$d^nF(0)$ is separately continuous on $\NC_\Ckl$ in its $n$ variables.
Hence by the nuclear theorem (e.g., [6], [13]) there exists $\Phi^{(n)}
\in(\NC_\Ckl^*)^{\sten n}$ so that
$$
\langle\Phi\nn,\xi_1\sten\cdots\sten\xi_n\rangle
        ={1\over n!}\,d^nF(0)(\xi_1,\ldots,\xi_n),\quad \xi_1,\ldots,\xi_n\in\NC_\Ckl.
$$
\def\pp{{p^\prime}}
\def\lang{\langle}
\def\rang{\rangle}
Let $\pp>p$ be such that the embedding $\iota:\,\HC_\pp\lra\HC_p$
is Hilbert--Schmidt, and let $(e_k,\,k\in\N)$ be an orthonormal basis
of $\HC_\pp$ in $\NC$. For $n\in\N$,  
$(e_{k_1}\ten\cdots\ten e_{k_n},\,k_i\in\N,\,i=1,\ldots,n)$ is an 
orthonormal basis of $\HC^{\ten n}_{\Ckl,\pp}$. Then we can estimate in
the following way (cf.\ (5)):
$$
\eqalign{
     |\Phi\nn|^2_{-\pp}&=\sum_{k_1,\ldots,k_n}
     |\lang\Phi\nn,e_{k_1}\ten\cdots\ten e_{k_n}\rang|^2\cr
  &=\sum_{k_1,\ldots,k_n}(n!)^{-2}\left|d^nF(0)(e_{k_1},\ldots,e_{k_n})\right|^2\cr
     &\le C^2(n!)^{-1}(2e^2K)^n(\sum_{k=1}^\infty|\iota e_k|_p^2)^n\cr
           &=C^2(n!)^{-1} (2e^2K\,\|\iota\|_{\rm HS}^2)^n\cr
           }
$$
i.e., $\Phi\nn\in\HC^{\sten n}_{\Ckl,-\pp}$, and
$$
\|\Phi\nn\|_{-\pp,-q}^2\le C^2 (2^{1-q}e^2K\,\|\iota\|_{\rm HS}^2)^n.
\eqno{(7)}
$$
For $\Phi$ given by the 
sequence $(\Phi\nn,\,n\in\N_0)$ ($\Phi^{(0)}\equiv F(0)$) we have
$$
\eqalign{
\|\Phi\|^2_{-\pp,-q}&\le C^2\sum_{n=0}^\infty(2^{1-q}e^2K\,\|\iota\|_{\rm HS}^2)^n\cr
        &=C^2(1-2^{1-q}e^2K\,\|\iota\|_{\rm HS}^2)^{-1}\cr
        &<+\infty,\cr
        }
$$
if we choose $q$ large enough so that $2^{1-q}e^2K\,\|\iota\|_{\rm HS}^2<1$.
In particular, $\Phi\in(\NC)^*$, and for $f\in\NC$ we have by (4),                  
$$
\eqalign{
S\Phi(f)&=\sum_{n=0}^\infty\lang\Phi\nn,f^{\ten n}\rang\cr
        &=\sum_{n=0}^\infty{1\over n!}\,\wh{d^nF(0)}(f)\cr
        &=F(f).\cr
        }
$$
Uniqueness of $\Phi=S^{-1}F$ follows from the fact that the exponential
vectors are total in $(\NC)$.\Kasten        
\medskip
As a by-product of the above proof we obtain the following localization
result for generalized functionals.
\def\r{\rho}
\Th{Corollary 12.}{\it
Given a U--functional $F$ satisfying C.2$\,^\prime$.
Let $\pp>p$ be such that the embedding $\iota:\,\HC_\pp\ra\HC_p$ is
Hilbert--Schmidt, and $q\in\N_0$ so that 
$\r:=2^{1-q}e^2K\,\|\iota\|_{\rm HS}^2<1$. 
Then $\Phi:=S^{-1}F\in\G_{-q}(\HC_{-\pp})$,
and
$$
\left\| \Phi \right\| _{-\pp,-q}\leq C(1-\r)^{-1/2}.
\eqno{(8)}
$$
}

For analogous results in white noise analysis see, e.g., \cite{23}, 
\cite{37}, \cite{50}.
\medskip
Within the framework established here one can treat the following and
numerous other examples in a unified way.
\mn
{\bf Example 13.}\quad
We choose the triplet%
$$
{\cal S}({\R}^n)\subset L^2({\R}^n)\subset {\cal S}^{\prime }({\R}^n), 
$$
and equip ${\cal S}^{\prime }({\R}^n)$ with the Gaussian
measure with characteristic functional 
$$
C(f)=e^{-\frac 12\int f^2(t)\,d^nt},\quad f\in {\cal S}\left( {\R}^n\right). 
$$
Then the framework allows to discuss functionals of white noise with $n$%
--dimensional time parameter \cite{48}.
\mn
{\bf Example 14.}\quad
If we choose a finite direct sum of identical copies of ${\cal S}(%
{\R})\,$ as the basic real nuclear space we obtain the characterization
of the space of Hida distributions of the
noise of an $n$-dimensional Brownian motion \cite{48}.
\medskip
We close this section by the corresponding characterization theorem for 
$\left( {\cal N}\right)$. This result is independently due to \cite{21}, 
\cite{31}, \cite{33}, and has been generalized and modified in various
ways, e.g., \cite{37}, \cite{50}, \cite{54}.

\Th{Theorem 15.}{\it
A mapping $F:{\cal N\rightarrow }{\C}$ is the $S$--transform of an
element in $({\cal N)}$ if and only if it admits C.1 and the following
condition
\medskip
\iitem{C.3} there exists a system of norms 
$( \left|\, \cdot\, \right|_{-p},\,p\in\N_0)$, which yields the
inductive limit topology on ${\cal N}^{*}$, and such that  
for all $p\geq 0$ and $\epsilon >0$ there exists 
$C_{p,\e}>0$ so that\par 
$$
\left| F(zf)\right| \leq C_{p,\e}\exp\left( \epsilon |z|^2
        \left|f\right| _{-p}^2\right),\quad f\in\NC,\,z\in\C.
\eqno{(9)}
$$ 
}
\medskip
If for $F$ conditions C.1 and C.3 are satisfied we say that $F$ is of 
order 2 and minimal type.
\mn{\sl Proof.}\quad If $\varphi \in ({\cal N)}$ then condition C.1 is 
satisfied as a consequence of Theorem 11. For any $p,q\geq 0$ we 
estimate as follows
$$
\eqalign{|S\vp(zf)|&=|\sum_{n=0}^\infty\lang\vp\nn,(zf)^{\ten n}\rang|\cr
        &\le \sum_{n=0}^\infty |z|^n|\vp\nn|_p\,|f|_{-p}^n\cr
        &\le (\sum_{n=0}^\infty n!\,2^{nq}\,|\vp\nn|_p^2)^{1/2}
         (\sum_{n=0}^\infty {1\over n!}\,(2^{-q}|z|^2|f|^2_{-p})^n)^{1/2}\cr
        &=\|\vp\|_{p,q}\,\exp(2^{1-q}|z|^2|f|_{-p}^2).\cr               
}
$$
Hence condition C.3, too, is necessary.
\medskip
Conversely, let $F$ be a $U$--functional of order 2 and minimal type.
From $F$, construct a sequence $\vp=(\vp\nn,\,n\in\N_0)$ of continuous
linear forms $\vp\nn$ on $\NC^{\sten n}$ as in the proof of Theorem 10.
We have to show that $\vp$ belongs to $\G_q(\HC_r)$ for all $r,q\in
\N_0$. Let $r,q\in\N_0$ be given. Choose $p>r$ such that the injection
$\iota:\,\HC_p\ra\HC_r$ is Hilbert--Schmidt. Then so is the injection
$\iota^*:\, \HC_{-r}\ra\HC_{-p}$. $\e>0$ in (9) is chosen so that
$\r:=\e2^{1+q}e^2\,\|\iota^*\|^2_{\rm HS}<1$. Then the analogue of (7) reads
$$
\eqalign{\|\vp\nn\|_{r,q}^2
        &\le C_{p,\e}^2 (2^{q+1}e^2\e\,\|\iota^*\|^2_{\rm HS})^n\cr
        &=C^2_{p,\e}\r^n,\cr}
$$
and we get
$$
\eqalign{\|\vp\|_{r,q}&=(\sum_{n=0}^\infty \|\vp\nn\|_{r,q}^2)^{\half}\cr
        &\le C_{p,\e}(1-\r)^{-\half}.\cr
        }
$$
Thus $\vp\in(\NC)$, and the proof is complete.\Kasten

\section{4. Corollaries}\noindent
One useful application of Theorem 11 is the discussion of
convergence of a sequence of generalized functionals. A first version of
this theorem is worked out in \cite{42}. Here we use our more general
setting to state

\Th{Theorem 16.}{\it
Let $(F_n,\,n\in\N)$ denote a sequence of $U$--functionals such that
\medskip
\iitem{1.} $(F_n(f),\,n\in\N)$ is a Cauchy sequence for all $f\in {\cal N}$,
\iitem{2.} there exists a continuous norm $\left|\, \cdot\, \right| $ 
   on ${\cal N}$ and $C,\, K>0$ such that 
   $\left| F_n(zf)\right| \leq Ce^{K\left|z|^2|f\right|^2}$ for all 
   $f\in {\cal N},\,z\in\C$, and for almost all $n\in\N$.
\mn 
Then $(S^{-1}F_n,\,n\in\N)$ converges strongly in $(\NC)^*$.
}

\mn{\sl Proof.}\quad The assumptions and inequality (8) imply
that there exist $p,q\geq 0$ and $\r\in(0,1)$ such that 
for all $n\in\N$,
$$
\left\| \Phi _n\right\| _{-p,-q}\leq C(1-\r)^{-\half} 
$$
where $\Phi _n=S^{-1}F_n$. Since ${\cal E}$ is total in 
$\Gamma _{-q}({\cal H}_{-p})$, assumption 1 implies that 
$(\lang\!\lang \Phi _n,\varphi \rang\!\rang,\,n\in\N)$ is a 
Cauchy sequence for all 
$\varphi \in ({\cal N})$. Since $\left( {\cal N}\right) ^{*}$ is the dual 
of the countably Hilbert space $(\NC)$, which is in particular Fr\'echet,
it follows from the Banach--Steinhaus theorem that $(\NC)^*$ is 
weakly sequentially complete. Thus there exists 
$\Phi \in \left({\cal N}\right) ^{*}\,$ such that $\Phi $ is the weak 
limit of $( \Phi_n,\,n\in\N)$. The proof is concluded by the remark
that weak and strong convergence of sequences coincide in the
duals of nuclear spaces (e.g., [13]).\Kasten
\medskip
As a second application we consider a theorem which concerns the integration
of a family of generalized functionals.
\def\L{\Lambda}
\Th{Theorem 17.}{\it
Let $\left( \L ,{\cal A},\nu \right) $ be a measure space, and $\lambda
\mapsto \Phi _\lambda $ a mapping from $\L $ to $\left( {\cal N}%
\right) ^{*}$. We assume that the $S$--transform $F_\lambda 
=S\Phi_\lambda $
satisfies the following conditions:
\medskip
\item{1.} for every $f\in{\cal N}$ the mapping $\lambda 
\mapsto F_\lambda \left( f\right)$ is measurable,
\item{2.} there exists a continuous norm $|\,\cdot\,|$ on $\NC$ so that
for all $\l\in\L,\,F_\lambda$ satisfies the bound $\left|F_\lambda (zf)\right| 
\leq C_\lambda e^{K_\lambda |z|^2|f|^2}$, and such that 
$\lambda \mapsto K_\lambda$ is bounded $\nu$--a.e., and $\lambda\mapsto
C_\lambda$ is integrable with respect to $\nu$. 
\mn
Then there are $q,p\geq 0
$ such that $\Phi_{\cdot}$ is Bochner integrable on 
$\Gamma _{-q}\left({\cal H}_{-p}\right) $. Thus in particular,
$$
\int_\L \Phi _\lambda\, d\nu (\lambda )\in \left( {\cal N}\right) ^{*}, 
$$
and%
$$
S\left( \int_\L \Phi _\lambda\, d\nu (\lambda )\right) (f)=\int_\L
S\Phi _\lambda (f)\,d\nu (\lambda ),\quad f\in\NC. 
$$
}

\mn{\sl Proof.}\quad In inequality (3) for $F_\l(zf)$ we
can replace $K_\l$ by its bound. With this modified estimate and 
Corollary 12
we can find $p,q\geq 0$ and $\r\in(0,1)$ such that for all $\l\in\L$, 
$$
\left\| \Phi _\lambda \right\| _{-p,-q}\leq C_\lambda(1-\r)^{-\half}.
\eqno{(10)}
$$
Since  the right hand side of (10) is integrable with respect to 
$\nu $, we only need to show the weak measurability of $\lambda \mapsto
\Phi _\lambda $ (see \cite{52}). But this is obvious because $\lambda
\mapsto \left\langle\! \left\langle \Phi _\lambda ,\varphi \right\rangle\!%
\right\rangle $ is measurable for all $\varphi \in {\cal E}$ which is total
in $\Gamma _q({\cal H}_p).$\Kasten
\bn
{\bf Acknowledgement.} We thank Professors S. W. He and
H. Sato for pointing out the gap in \cite{42}, and Professor B. \O ksendal
for helpful discussions. We owe special thanks to Professor L.I. Ronkin
who taught us about Bernstein's theorem and its generalizations. 
\goodbreak        
\input refmac
\def\bibitem#1{\li{[#1]}}
\def\li{\par\noindent
                      \hangindent=1\parindent
                      \ltextindent}
\li{[1]} Akhiezer, N.I. and Ronkin, L.I. (1973), {\it On separately
analytic functions of several variables and theorems on ``the thin
edge of the wedge"}, Russian Math.\ Surveys 28, No.\ 3, 27--44
\bibitem{2}  Albeverio, S., Hida, T., Potthoff, J. and Streit, L. (1989),%
{\it \ The vacuum of the Hoegh-Krohn model as a generalized white noise
functional,} Phys. Lett. B 217, 511--514.

\bibitem{3}  Albeverio, S., Hida, T., Potthoff, J., R\"ockner, M. and
Streit, L. (1990){\it , Dirichlet forms in terms of white noise analysis I:
Construction and QFT examples.} Rev. Math. Phys. 1, 291--312

\bibitem{4}  Albeverio, S., Hida, T., Potthoff, J., R\"ockner, M. and
Streit, L. (1990),{\it \ Dirichlet forms in terms of white noise analysis
II: Closability and diffusion processes}. Rev. Math. Phys. 1, 313--323.

\bibitem{5}  Albeverio, S., Kondratiev, Yu.G. and Streit, L. (1993), {\it %
How to generalize White Noise Analysis to Non--Gaussian Spaces}. In: Dynamics
of Complex and Irregular Systems. Ph. Blanchard et al., eds. World
Scientific.

\bibitem{6}  Berezansky, Yu. M. and Kondratiev, Yu. G. (1988), {\it %
Spectral Methods in Infinite--Dimen\-sio\-nal Analysis}, (in Russian), Naukova
Dumka, Kiev.

\bibitem{7}  Cochran, G., Lee, J.--S. and Potthoff, J. (1993): {\it %
Stochastic Volterra equations with singular kernels}; Preprint.

\bibitem{8}  Colombeau, J.-F. (1982), {\it Differential calculus and
holomorphy. }Mathematical Studies 64, North--Holland, Amsterdam.

\bibitem{9}  Cook, J. (1953), {\it The mathematics of second quantization.%
} Trans. Amer. Math. Soc. 74, 222--245.

\bibitem{10}  de Faria, M., Hida, T., Streit, L., and Watanabe, H., {\it %
Intersection local times as Generalized White Noise Functionals}, in
preparation.

\bibitem{11}  de Faria, M., Potthoff, J. and Streit, L. (1991), {\it The
Feynman integrand as a Hida distribution.} J. Math. Phys. 32, 2123--2127.

\bibitem{12}  Dineen, S. (1981), {\it Complex Analysis in Locally Convex
Spaces,} Mathematical Studies 57, North--Holland, Amsterdam.

\bibitem{13}  Gel'fand, I.M. and Vilenkin, N.Ya. (1968),{\it \ Generalized
Functions}, Vol. IV, Academic Press, New York and London.

\bibitem{14}  Hida, T. (1989), {\it Infinite-dimensional rotation group
and unitary group.} Lecture Notes in Math 1379, 125--134.

\bibitem{15}  Hida, T., Potthoff, J. and Streit, L.(1988), {\it Dirichlet
Forms and white noise analysis.} Commun. Math. Phys. 116, 235--245.

\bibitem{16}  Hida, T. and Streit, L. (1983), {\it Generalized Brownian
functionals and the Feynman integral}. Stoch. Proc. Appl. 16, 55--69.

\bibitem{17}  Hida, T., Kuo, H.--H., Potthoff, J. and Streit, L. (1993),{\it %
\ White noise. An infinite dimensional calculus}. Kluwer, Dordrecht.

\bibitem{18}  Holden, H., Lindstr\o m, T., \O ksendal, B., Ub\o e, J. and
Zhang, T.--S. (1993), {\it Stochastic boundary value problems: A white noise
functional approach}; Probab. Th. Rel. Fields { 95}, 391--419.

\bibitem{19}  Holden, H., Lindstr\o m, T., \O ksendal, B., Ub\o e, J. and
Zhang, T.--S.: {\it ,} Preprint (1993).

\bibitem{20}  Khandekar, D.C. and Streit, L. (1992),{\it \ Constructing the
Feynman integrand}. Ann. Physik 1, 49--55.

\bibitem{21}  Kondratiev, Yu.G. (1980), {\it Nuclear spaces of entire
functions of an infinite number of variables, connected with the rigging
of a Fock space.} In: Spectral Analysis of Differential Operators.
Math.\ Inst., Acad.\ Sci.\ Ukrainian SSR, p.\ 18--37. English translation:
Selecta Math.\ Sovietica 10 (1991), 165--180.

\bibitem{22}  Kondratiev, Yu.G. and Samoilenko, Yu.S. (1976),{\it \
Integral representation of generalized positive definite kernels of an
infinite number of variables}. Soviet Math. Dokl. 17, 517--521.

\bibitem{23}  Kondratiev, Yu.G. and Streit, L. (1991a),{\it \ A remark
about a norm estimate for White Noise distributions.} Ukrainian Math. J.

\bibitem{24}  Kondratiev, Yu.G. and Streit, L. (1992), {\it Spaces of White
Noise distributions: Constructions, Descriptions, Applications.} I. BiBoS
preprint no. 510, to appear in Rep.Math.Phys.

\bibitem{25}  Kr\'ee, P. (1988), {\it La theorie des distributions en
dimension quelconque et l~`integration stochastique}. Lecture Notes in Math.
1316, 170--233, Springer, Berlin, Heidelberg, New York.

\bibitem{26}  Kubo, I. and Takenaka, S. (1980a), {\it Calculus on Gaussian
white noise I.\/} Proc. Japan Acad. 56, 376--380.

\bibitem{27}  Kubo, I. and Takenaka, S. (1980b), {\it Calculus on Gaussian
white noise II.\/} Proc. Japan Acad. 56, 411--416.

\bibitem{28}  Kubo, I. and Takenaka, S. (1981), {\it Calculus on Gaussian
white noise III.\/} Proc. Japan Acad. 57, 433--437.

\bibitem{29}  Kubo, I. and Takenaka, S. (1982), {\it Calculus on Gaussian
white noise IV.\/} Proc. Japan Acad. 58, 186--189.

\bibitem{30}  Kuo, H.--H. and Potthoff, J.: {\it Anticipating stochastic
integrals and stochastic differential equations;} in White Noise Analysis --
Mathematics and Applications, eds.: T. Hida, H.--H. Kuo, J. Potthoff and L.
Streit. Singapore: World Scientific (1990).

\bibitem{31}  Kuo, H.-H., Potthoff, J. and Streit, L. (1991), {\it A
characterization of white noise test functionals.} Nagoya Math. J. {\bf 121}%
,185--194.

\bibitem{32}  Lascheck, A., Leukert, P., Streit, L., Westerkamp, W.
(1993) {\it %
Quantum mechanical propagators in terms of Hida distributions}. 
Rep.\ Math.\ Phys.\ 33, 221--232

\bibitem{33}  Lee, Y.J. (1989), {\it Generalized Functions of Infinite
Dimensional spaces and its Application to White Noise Calculus}. J. Funct.
Anal. {\bf 82}, 429--464.

\bibitem{34}  Meyer, P.A. and Yan, J.-A. (1990),{\it \ Les ``fonctions
caract\'eristiques'' des distributions sur l'\'espace de Wiener.} Seminaire de
Probabilites XXV, ed.: J Azema, P.A. Meyer, M. Yor, Springer, p. 61--78.

\bibitem{35}  Nachbin, L. (1969), {\it Topology on spaces of holomorphic
mappings.} Springer, Berlin.

\bibitem{36}  Nelson, E. (1973), {\it Probability theory and Euclidean
quantum field theory. In: Constructive Quantum Field Theory,} ed. by G. Velo
and A. Wightman, Springer, Berlin, Heidelberg, New York.

\bibitem{37}  Obata, N. (1991), {\it An analytic characterization of
symbols of operators on white noise functionals. }Journal of Functional
Analysis.

\li{[38]} Pietsch, A. (1969), {\it Nukleare Lokal Konvexe R\"aume},
  Berlin, Akademie Verlag.

\bibitem{39}  Potthoff, J. (1991), {\it Introduction to white noise
analysis}; in Control Theory, Stochastic Analysis and Applications, S. Chen,
J. Yong (ed.s), Singapore, World Scientific.

\bibitem{40}  Potthoff, J. (1992), {\it White noise methods for stochastic
partial differential equations.} In: Stochastic Partial Differential
Equations and Their Applications, B.L. Rozovskii, R.B. Sowers (ed.s),
Berlin, Heidelberg, New York, Springer.

\bibitem{41}  Potthoff, J.: {\it White noise approach to parabolic
stochastic differential equations}, in preparation.

\bibitem{42}  Potthoff, J. and Streit, L. (1991), {\it A characterization of
Hida distributions.} J. Funct. Anal. {\bf 101}, 212--229.

\bibitem{43}  Potthoff, J. and Streit, L. (1993), {\it Invariant states on
random and quantum fields: }$\phi -${\it bounds and white noise analysis}.
J. Funct. Anal. {\bf 101}, 295--311.

\li{[44]} Reed, M., and Simon, B. (1972), {\it Methods of Modern Mathematical 
Physics I: Functional Analysis}. Academic Press,  New York and London.

\li{[45]} Siciak, J. (1969), {\it Separately analytic functions and
envelopes of holomorphy of some lower dimensional subsets of $\C^n$},
Ann.\ Polonici Math.\ 22, 145--171

\bibitem{46}  Simon, B. (1974), {\it The P}$\left( \phi \right) _2${\it \
Euclidean (quantum) field theory}. Princeton University Press, Princeton.

\bibitem{47}  Segal, I. (1956), {\it Tensor algebras over Hilbert spaces}.
Trans. Amer. Math. Soc. {\bf 81}, 106--134.

\bibitem{48}  Streit, L. and Westerkamp, W. (1993),{\it \ A generalization
of the characterization theorem for generalized functionals of White Noise}.
In: Dynamics of Complex and Irregular Systems. Ph. Blanchard et al., eds.
World Scientific.

\bibitem{49}  Watanabe, H. (1991), {\it The local time of self-intersections
of Brownian Motions as generalized Brownian functionals}, Lett. Math. Phys. 
{\bf 23}, 1--9.

\bibitem{50}  Yan, J.-A. (1990), {\it A characterization of white noise
functionals.} Preprint.

\bibitem{51}  Yokoi, Y. (1993), {\it Simple setting for white noise
calculus using Bargmann space and Gauss transform.} Preprint.

\bibitem{52}  Yosida, K. (1980), {\it Functional Analysis}. Springer,
Berlin.

\li{[53]} Zaharjuta, V.P. (1976), {\it Separately analytic functions,
generalizations of Hartogs' theorem, and envelopes of holomorphy},
Math.\ USSR Sbornik 30, 51--67

\bibitem{54}  Zhang, T.--S.: {\it Characterization of white noise test
functions and Hida distributions.} Stochastics {\bf 41} (1992) 71--87.

\end

%% file: mydef.tex

\magnification = \magstep1 
\font\small=cmr10 at 10truept

\font\large=cmbx10 scaled \magstep2
\font\spec=cmcsc10 scaled \magstep1
\font\bi=cmbxti10 scaled \magstep1

\def\today{\ifcase\month\or January\or February\or March\or
April\or May\or June\or July\or August\or September\or
October\or November\or December\fi\space\number\day,\  
\number\year}

\def\heute{\number\day.\ \ifcase\month\or Januar\or Februar\or M\"arz\or
April\or Mai\or Juni\or Juli\or August\or September\or
Oktober\or November\or Dezember\fi\space\number\year}

\overfullrule=0pt
\def\leaderfill{\leaders\hbox to 1em{\hss.\hss}\hfill}
\def\part#1{\bigskip\bigskip\bigskip\noindent {\large #1}
            \bigskip\bigskip}

\def\section#1{\bigskip\medskip\goodbreak\noindent {\bf #1} \bigskip}

\def\Th#1{\bigskip\goodbreak\noindent {\bf #1}\quad }

\def\Kasten{\hfill\hbox{             
    \vrule width 0.3pt height 4.7pt depth 0.3pt
    \vrule width 4.4pt height 0pt   depth 0.3pt \hskip -4.4pt
    \vrule width 4.4pt height 4.7pt depth -4.4pt
    \vrule width 0.3pt height 4.7pt depth 0.3pt \hskip 2pt}}
\def\1{{1\!\!1}}
\def\3{\ss }

\def\bn{\bigskip\noindent}
\def\C{\hbox{ \vrule width 0.6pt height 6.0pt depth 0pt \hskip -3.5pt}C}
\def\Ckl{{\,\hbox{\vrule width 0.4pt height 4.3pt depth 0pt
               \hskip -2.5pt}C}}
\def\Ckkl{{\,\hbox{\vrule width 0.3pt height 3.2pt depth 0pt
               \hskip -2.2pt}C}}

\def\cline#1{\centerline{#1}}

\def\EC{{\cal E}}

\def\e{\varepsilon }

\def\G{\Gamma }

\def\HC{{\cal H}}
\def\half{{1 \over 2}}

\def\iitem{\itemitem}

\def\KC{{\cal K}}

\def\L{\Lambda }

\def\l{\lambda }

\def\lms{\longmapsto }
\def\lra{\mathop{\longrightarrow}}

\def\mn{\medskip\noindent}

\def\N{{I\!\! N}}
\def\NC{{\cal N}}

\def\nn{^{(n)}}

\def\R{{I\!\!R}}

\def\ra{\rightarrow }

\def\SC{{\cal S}}

\def\sn{\smallskip\noindent}
\def\sr{{\SC(\R)}}
\def\ssr{{\SC'(\R)}}

\def\ten{\otimes}
\def\sten{{\widehat\otimes}}

\def\vp{\varphi }

\def\wh{\widehat}

\def\Z{Z\!\!\! Z }

%% file: refmac.tex
\par
\narrower
\def\li{\par\noindent
                      \hangindent=3\parindent
                      \ltextindent}
\def\ltextindent#1{\hbox to \hangindent{#1\hss}\ignorespaces}

\def\item{\li}

\def\iitem{\li}

\font\small=cmr10 at 10truept
\font\it =cmti10 at 10 truept
\font\bf=cmbx10 at 10 truept

\bigskip\bigskip\noindent
{\bf References.}
\medskip\noindent
\baselineskip 10pt
\small